\documentclass[12pt,amssymb]{article}
\usepackage{amsmath,amssymb}
\usepackage{amscd}
\newtheorem{Theorem}{Theorem}
\newtheorem{Theorem2}{Theorem}[section]
\newtheorem{Lemma}[Theorem2]{Lemma}
\newtheorem{Proposition}[Theorem2]{Proposition}
\newtheorem{Corollary}[Theorem2]{Corollary}
\newtheorem{Remark}[Theorem2]{Remark}
\newtheorem{Definition}[Theorem2]{Definition}
\newtheorem{Notation}[Theorem2]{Notation}
\newcommand{\pf}{{\bf Proof}\ }

\newcommand{\Q}{{\mathbb Q}}
\newcommand{\Z}{{\mathbb Z}}
\newcommand{\N}{{\mathbb N}}
\newcommand{\C}{{\mathbb C}}
\newcommand{\proZ}{\widehat{\Z}}
\newcommand{\rem}{\begin{Remark}\rm }
\newcommand{\erem}{\end{Remark}}

\newcommand{\HH}{{\bf H}}

\newcommand{\nt}{\begin{Notation}\rm }
\newcommand{\ent}{\end{Notation}}

\newcommand{\Ker}{{\rm ker}\ }
\newcommand{\Dom}{{\rm Dom\ }}

\newcommand{\pr}{{\rm pr}}
\newcommand{\ex}{{\rm ex}}
\newcommand{\tp}{{\rm tp}}

\newcommand{\acl}{{\rm acl}}
\newcommand{\cl}{{\rm cl}}

\newcommand{\ssn}{\section}
\newcommand{\df}{\begin{Definition}\rm }
\newcommand{\edf}{\end{Definition}}
\newcommand{\bl}{\begin{Lemma}}
\newcommand{\el}{\end{Lemma}}
\newcommand{\bt}{\begin{Theorem}}
\newcommand{\et}{\end{Theorem}}
\newcommand{\bp}{\begin{Proposition}}
\newcommand{\ep}{\end{Proposition}}
\newcommand{\be}{\begin{equation}}
\newcommand{\ee}{\end{equation}}
\newcommand{\inv}{^{-1}}
\newcommand{\bc}{\begin{Corollary}}
\newcommand{\ec}{\end{Corollary}}
\newcommand{\lb}{\label}

\newcommand{\ra}{\right\rangle}
\newcommand{\la}{\left\langle}
\newcommand{\subs}{\subseteq}

\newcommand{\qed}{$\Box$ \\ \\}
\newcommand{\smin}{\setminus}

\newcommand{\To}{\longrightarrow}
\newcommand{\ov}{\overset}
\begin{document}
\title{Covers of the multiplicative group of an algebraically closed
field of characteristic zero}
\author{Boris Zilber\\ Unversity of Oxford}
\date{October 21, 2003}
\maketitle
%\footnotetext{{\em Mathematical Subject Classification} 12F10, 12L12, 03C60}
\ssn{ Introduction and results}

Consider the classical universal cover of the one dimensional
complex torus $\C^*,$ which gives us the exact sequence
\be \lb{e04} 0\To \Z\ov{i}{\To} \C\ov{\exp}{\To} \C^*\to 1. \ee

Model theoretically one can interprete the sequence  as a structure in various
ways. The simplest 
 algebraic structure on the sequence which  bears an 
interesting algebro-geometric information
is the one with the additive group structure in the middle and with the full algebraic geometry on $\C^*.$
The latter is equivalent to treating $\C^*$ as $\C\smin \{ 0\}$ with the full field structure on $\C.$ We call this structure {\bf a group cover} of the multiplicative group of the field. 
   
Is  the group cover of $\C^*$ determined uniquely? We can put the
question in the following precise form:

{\em Given an exact sequence
\be \lb{e05} 0\To \Z\overset{i_H}{\To} H\ov{\ex}{\To} \C^*\to 1 \ee
with $H$ a torsion-free divisible abelian group and $\ex$ a group
homomorphism, is there an isomorphism $\sigma$ between
the groups covers (\ref{e04}) and (\ref{e05}), that is a group-isomorphism $\sigma_H:\C\to H$  and
a field automorphism $\sigma_C: \C\to \C$ such that  
 all the arrows commute? }

$$
\begin{CD}
0  @>>> \Z @>{i_C}>> \C     @>{\exp}>>   \C^* @>>> 1 \\
@.      @VViV       @VV\sigma_HV             @VV\sigma_CV       @.\\
0  @>>> \Z @>{i_H}>> H  @>{\ex}>> \C^* @>>> 1
\end{CD}
$$

This paper  gives a positive answer to the question in a more general form.

\bt \lb{c1} Let $F^{\times}$ be the multiplicative group of an algebraically closed field of
characteristic 0, $H$ an abelian divisible torsion free group such
that the sequence
\be \lb{e4} 0\To \Z\ov{i}{\To} H\ov{\ex}{\To} F^{\times}\To 1 \ee
is exact.  Then the isomorphism type of the sequence is determined by
the isomorphism type of the field 
$F.$ In other words, if 
\be \lb{e5} 0\To \Z\ov{i'}{\To} H'\ov{\ex'}{\To} F'^{\times}\To 1\ee   
is another such sequence, with a field $F'$ isomorphic to $F,$  then there is an isomorphism $\sigma$ between
the sequences, which is a field isomorphism on $F,$ a group
isomorphism on $H$ and $\Z,$ and all the arrows commute. 
\et

It is worth reminding here that the isomorphism type of an algebraically
closed field of a given characteristic is, by Steinitz' theorem,  determined by its
transcendence degree. \\

The results of the paper furnish an algebraic  background for the
study  of a more complicated structure, {\em a field with
pseudo-exponentiation}  (see \  [Z1, section 2] and [Z3]) which takes into account
that the cover $H$ bears field structure as well. 
 
The proof of Theorem~\ref{c1} is a model-theoretic consequence of   
Theorem~\ref{t0} (see the introduction below)  formulated and proved in 
a field-theoretic  language.  The formulation is rather
 technical but, importantly, it is an arithmetic
statement equivalent to the geometric form of Theorem~\ref{c1}. The
equivalence is due to the model-theoretical Keisler-Shelah  theory of 
{\em excellency.} We do not give full details for this equivalence here leaving
it to the forthcoming paper [Z4].

Theorem~\ref{c1} is not trivial due to the fact that, once for $a\in
F^{\times}$ we fix an $h\in H$ such that $ex(h)=a,$ we fix the whole subgroup
$\ex(\Q\cdot h).$  In particular if, for
example, $F=F'$ we can not in general take $\sigma$ to be  identity
on $F.$ 

On the other hand, there is a version of the statement which is quite
easy to prove; this is the case when in the definitions  $\Z$ is replaced 
by  $\hat \Z,$ the completion of $\Z$ in the profinite topology.
This corresponds to the well-known SGA-construction of ``the algebraic
$\pi_1$''. We discuss these issues and provide 
detailed proofs  in the last section of the
paper.   \\ \\

The author would like to thank Neil Dummigan for a very substantial
contribution to this paper. On the initial step of the research he 
directed the author to  certain 
Galois calculations, which developed later into the proof of the
crucial Lemma~\ref{l3}, and also took  part in the proof of an important
Lemma~\ref{l0}, not mentioning helpful general discussions and
suggestions. 

The author is very thankful to D.Pierce for carefully reading the
paper and suggesting many valuable improvements.
Also thanks are due to P.Voloch who provided us with a proof of a
version of Corollary~\ref{c0} over an algebraically closed field. 
\\ \\

In the rest of the section we introduce our main definitions and discuss the main results on a more technical level.\\   

\nt For a non-torsion $a\in \C^{\times}$ let $a^{\Q}$ be a (non unique)
multiplicative subgroup of $\C^{\times}$ containing $a$ and isomorphic to the
additive group $\Q$ of the rational numbers. We call such a subgroup
{\bf a multiplicatively divisible subgroup associated with $a$}. (An
example of such a subgroup one gets by fixing a value for $\ln a$ and
letting $$a^{\Q}=\{ \exp(q\ln a): q\in \Q\}.)$$ 
Notice that the notation $a^{\Q}$ makes also sense when $a$ is a root
of unity, then this is just the torsion group of the field.

We further denote\\

 $\la a_1,\dots, a_n\ra$ is the multiplicative subgroup generated by
 $a_1,\dots, a_n\in \C^{\times}.$\\

$\mu_n$
 is the group
of  roots of unity of power $n$ in $\C^{\times},$\\

$\mu$
 is the group
of all the roots of unity in $\C^{\times}.$\ent

\rem  Notice that, for a fixed $a\in \C^{\times}$ as above, a choice of $a^{\Q}$ can be
done by choosing  roots $a^\frac{1}{m}$ of $a$ of powers $m\in
\N$ in an agreeing way, which thus form a projective system, in a
natural bijection to the projective system
$$\mu_{ml}\to^{x^m}\mu_l$$ of subgroups of roots of unity. The
projective limit of the
latter can be identified with the group $\hat \Z,$ the completion of the
cyclic group $\Z$ in the profinite topology. Thus any choice of
$a^{\Q}$ corresponds to a point in $\hat \Z.$\erem

When $a^{\Q}$ is given, $a^\frac{1}{m}$ will stand for the unique root of $a$ of
power $m$ in $a^{\Q}.$

For $P$ a subfield of $\C$ and $X_1,\dots,X_n\subs \C$ we let
$P(X_1,\dots,X_n)$ be the subfield generated by $X_1,\dots,X_n$ over
$P.$ \\

We say that $\{ b_1,\dots, b_l\}\subs \C^{\times}$ 
{\bf multiplicatively independent} if, for  integer $n_1,\dots, n_l$ 
$$b_1^{n_1}\cdot\dots \cdot b_l^{n_l}=1\mbox{  if and only if   }n_1=\dots=n_l=0.$$

\df  Given some subgroups  $b_1^{\Q},\dots, b_l^{\Q}$
 of $\C^{\times}$ a subfield $F\subs \C$ and some positive integer $m,$
we will say that the group elements
{\bf  $b_1^{\frac{1}{m}}\in b_1^{\Q},\dots,
 b_l^{\frac{1}{m}}\in b_l^{\Q}$ 
 determine the isomorphism type of $b_1^{\Q},\dots\cdot, b_l^{\Q}$
over $F$} if, for any subgroups $c_1^{\Q},\dots, c_l^{\Q}$ of $\C^{\times},$
 any
field isomorphism 
$$\phi_m: F( b_1^{\frac{1}{m}},\dots, b_l^{\frac{1}{m}})\to F( c_1^{\frac{1}{m}},\dots, c_l^{\frac{1}{m}}),$$
$$b_i^{\frac{1}{m}}\mapsto c_i^{\frac{1}{m}}$$
over $F$   extends to 
a field isomorphism 
$$\phi_{\infty}:F(b_1^{\Q},\dots, b_l^{\Q}) \to F(c_1^{\Q},\dots, c_l^{\Q})$$
 $$b_i^{\Q}\mapsto c_i^{\Q}.$$\edf

%Notice that $\phi_{\infty}$ in the definition is, in general, not  unique .

\bt  \lb{t0} 
%Let $L\subs \C$ be an algebraically closed subfield or $L=\Q.$ 
Let
$P\subs \C$ be a finitely generated extension of $\Q$ and $L_1,\dots,
L_n$ algebraically closed subfields of the algebraic closure of $P,$
$n\ge 0.$
Let $a_1,\dots, a_r\in P^{\times},$    
   $b_1,\dots, b_l\in \C^{\times}$ and
$a_1^{\Q},\dots,a_r^{\Q},b_1^{\Q},\dots, b_l^{\Q}$ are multiplicatively divisible subgroups
 associated with the elements.
Suppose that 
 $b_1,\dots b_l$ are
 multiplicatively independent over  $\la a_1,\dots,
 a_r\ra \cdot L_1^{\times}\cdot\dots \cdot L_n^{\times}.$ 

    Then
 there is an $m\in \N$ such that $b_1^{\frac{1}{m}},\dots,
 b_l^{\frac{1}{m}}$ 
 determine the isomorphism type of $b_1^{\Q},\dots,b_l^{\Q}$
 over the field $\hat P=P(L_1,\dots,L_n,\mu,a_1^{\Q},\dots,
 a_r^{\Q}).$   
\et 

\rem  $\mu$ is needed in the definition of $\hat P$ only in case $n=0.$\erem

Theorem~\ref{t0} has the following
\bc \lb{c0} Let $W^{\frac{1}{n}}$ be the locus of (the minimal algebraic
variety containing) $\la b_{1}^{\frac{1}{n}},\dots,b_l^{\frac{1}{n}}\ra$ 
over $\hat P.$ Then 
for some $m$, for each $d$, the algebraic set
  \begin{equation*}
      \{(x_1,\dots,x_{\ell}):(x_1^d,\dots, x_{\ell}^d)\in W^{1/m}\}
  \end{equation*}
is irreducible over $\hat P$ and is precisely $W^{1/dm}$.
\ec

The proof of Theorem~\ref{c1} in  case $F$ is countable follows
directly from Theorem~\ref{t0}, with $n=0,$ by the standard back-and-forth
construction of the isomorphism. For this we only need to notice that if a
partial linear isomorphism $\sigma_H: H\to H'$ maps a $\Q$-subspace generated by linearly independent
$h_1,\dots,h_r\in H$ to the subspace generated by linearly independent $h'_1,\dots,h'_r\in H',$ with $a_1=\ex(h_1),\dots, a_r=\ex(h_r)$ and
$\sigma_F:a_i\to a'_i=\ex(h'_i)$ a Galois isomorphism, then we can extend  
$\sigma$ to any new $h.$ Indeed, let $\ex(h)=b,$ and choose  $m$ as
in Theorem~\ref{t0} and 
${b}^{\frac{1}{m}}=\ex(\frac{h}{m}).$ Then extend the field-isomorphism
$\sigma_F$ by defining ${b'}^{\frac{1}{m}}=\sigma_F({b}^{\frac{1}{m}})$
and $h'\in H'$ so that ${b'}^{\frac{1}{m}}=\ex(\frac{h'}{m}),$
and finally put $\sigma_H(h)=h'.$

The case of cardinality $\aleph_1$ can be done applying the case $n=1$ of  Theorem~\ref{t0} by using one more standard
model-theoretic trick. As we want to
 extend an isomorphism between two countable
sequences (\ref{e4}) one transcendence degree up $\aleph_1$ times,
we need to be able to extend $\sigma_0=(\sigma_{H_0},\sigma_{F_0})$ for countable $H_0$ and $F_0$ 
as above to $(H, F)$ with $F=\ex(H)$ an algebraically closed field of transcendence degree one over $F_0,$ say $$F=\acl F_0(b_0)= F_0(b_0,\dots, b_i,\dots).$$ 
We extend $\sigma_0$ by induction  to $\sigma_i=(\sigma_{H_i},\sigma_{F_i})$ 
with $F_i$ containing 
$b_0,\dots, b_{i-1}$
 and $H_i$ containing $h_0,\dots, h_{i-1},$ such that $\ex(h_0)=b_0,\dots,\ex(h_{i-1})=b_{i-1}.$  Since $b_0$ is transcendental over $F_0$ any choice of $b'_0\in F'$ and 
$h'_0\in H'$ with $\ex(h'_0)=b'_0$ will do for $\sigma_{H_1}(h_0)$ and 
$\sigma_{F_1}(b_0).$  

For  step $i>1$ let $L$ be the algebraically closed subfield of $F_0$ of finite transcendence degree of the form $F_0\cap \acl(b_0,\dots,b_i).$ 
Then, since $F_0$ and $\acl(b_0,\dots,b_i)$ are linearly disjoint over $L$  (see [L], Ch.VIII)
any Galois automorphism of $L(b_0^{\Q},\dots,b_i^{\Q})$ over $L$ can be extended to a  Galois automorphism of $F_0(b_0^{\Q},\dots,b_i^{\Q})$ over $F_0.$
In particular, if the isomorphism type of $b_0^{\Q},\dots,b_i^{\Q}$ over
$L$ is determined by $b_0^{\frac{1}{m}},\dots,b_i^{\frac{1}{m}}$ then 
the same is true over $F_0.$ Since such an $m$ is given by case $n=1$ of Theorem~\ref{t0}, we can proceed as above extending the isomorphism.

Surprisingly, this does not easily generalise to arbitrary
cardinalities. So, the passage from the countable to the general case
goes via a direct application of the main result of [Z2], a special
case of Shelah's {\em excellency }  theory [Sh], which gives
very general conditions for an isomorphism between uncountable 
structures to exist. 

It would be useful to remark that the linear disjointness argument used above corresponds to
the general notion of {\em splitting}, a part of the theory of excellency.
\\

It is highly desirable in view of the discussion of other
pseudo-analytic structures in [Z1]: 

- to generalise Theorem~\ref{t0} to fields of arbitrary characteristics;

- to prove a version of Theorem~\ref{c1} for the sequences of the form

\be \lb{e6} 0\To \Ker\ov{i}{\To} H\ov{\ex}{\To} A(F)\To 0, \ee
 
where $A(F)$ is the group of $F$-points of an abelian variety $A$ of
dimension $d$ with the field of definition $F_0\subs F,$ \ $F$ an 
algebraically closed field, \ $H$ an abelian torsion free group,
$\Ker$ an abelian group of rank $2d,$ and the isomorphism $\sigma$
on $A(F)$ is induced by an isomorphism of the field $F$ fixing $F_0.$

\ssn{ Proof of the main theorem }

Notice first that the above definition can be equivalently given as follows:

    { \em $b_1^{\frac{1}{m}}\in b_1^{\Q},\dots,
 b_l^{\frac{1}{m}}\in b_l^{\Q}$ 
 determine the isomorphism type of $b_1^{\Q}\cdot\dots\cdot b_l^{\Q}$
    over $F$ iff,}
 {\em given subgroups of the form $c_1^{\Q}\cdot\dots\cdot c_l^{\Q}$
and  a
field isomorphism $\phi_m$ over $F$ such that 
$\phi_m(b_i^{\frac{1}{m}})=c_i^{\frac{1}{m}},$ we can for any $d\in \N$
extend
 $\phi_m$  to 
a field isomorphism 
$$\phi_{dm}:F(b_1^{\frac{1}{dm}}, \dots, b_l^{\frac{1}{dm}})
 \to F(c_1^{\frac{1}{dm}},\dots, c_l^{\frac{1}{dm}})$$ taking
 $b_i^{\frac{1}{dm}}$ to $c_i^{\frac{1}{dm}}$}.\\

Indeed, though the latter formulation is formally weaker, it allows
us to define $$\phi_{\infty}=\bigcup_{d\in \N}\phi_{dm}.$$

 We prove Theorem~\ref{t0} through the following series of
 lemmas. Without loss of generality we assume that  $L_j\cap P$
 contains a transcendence basis of $L_j$ for each $j=1,\dots,n,$ that
 is $L_j$ is algebraic over its subfield $L_j\cap P.$  

\bl \lb{l0} The multiplicative group of the  field $P$ is
isomorphic to a direct product $A\cdot C$  of a free abelian group $A$
and  the group $C=P^{\times}\cap \mu.$ 

The multiplicative group of the  field $P(L),$ for $L$ an
algebraically closed field,  is
isomorphic to a direct product $A\cdot C$  of a free abelian group $A$
and the group  $C=L^{\times}.$    
\el
\pf Let first $P$ be a finite (algebraic) extension of $\Q.$ Then by
the theory of fractional ideals (see [He]), letting  $U$ be  the group
of units of the ring of integers of $P,$ the quotient-group
$ P^{\times}/U$
is a subgroup of the group of fractional ideals, which is a free
abelian group generated by the prime ideals. Hence $P^{\times}/U$ is free
too. 
By Dirichlet's Unit Theorem, $U=C\cdot U'$ where $C$ is the torsion subgroup of
$U$ and $U'$ a finitely generated free abelian group. 
 Hence $P^{\times}/C$ is free. It follows that
$P^{\times}=A\cdot C$ for some free $A.$

The multiplicative group of $P(L),$ for $P$ a finitely
generated extension of $\Q,$ is the multiplicative group of a field
which is of finite transcendence
degree over the algebraically closed field $L.$ Such a field  can be
viewed as the field of rational functions
of an algebraic variety over $L.$
According
to the normalization theorem ([H], I, ex.3.17), we can assume
 the variety normal. On a normal variety the concept of Weil divisor
([H], II.6)
makes sense, and the
divisors of functions (the principal divisors) form a subgroup of the 
free abelian group of all
Weil divisors on the variety. So $P(L)^{\times}/U,$ where $U$ is the subgroup of
 functions with trivial divisors, which is the group $L^{\times}$ of constant 
functions, is isomorphic to a subgroup $A$ of a free abelian
 group, and thus is free itself.

In the case that  
 $P$ is a finitely generated extension of $\Q$ we can 
 view $P$ as the function field
of an algebraic variety, over some finite extension of $\Q.$ 
We can again assume
 the variety normal over the algebraic closure of $\Q,$ by the
 normalisation theorem.
As the normalization procedure uses only
finitely many coefficients, so rather than going to the algebraic closure of
$\Q,$ some sufficiently large finite extension will do. So let the new $P$
be the function field of the variety over this sufficiently large finite
extension $K$ of $\Q.$
As above the group of principal Weil divisors $P^{\times}/K^{\times},$ is a free abelian
 group, say $A'.$
But $K^{\times},$   the multiplicative group
of a finite extension of $\Q,$  is isomorphic to $A''\cdot C,$ $A''$ free, by the above
proved. By the theory of abelian groups 
we thus
 have  $P^{\times}= A'\cdot A''\cdot C,$ with $A=A'\cdot A''$ free and the
 products direct.
\qed
 
 Recall that a subgroup $B$ of an abelian group is called {\bf
pure} if whenever the equation $x^n=b$  for a $b\in B$ and $n\ge 1$
has a solution in $A,$ it has a solution 
 in $B.$

\df A tuple $\{ a_1, \dots,a_k\}$ of elements of $P^{\times}$
will be called {\bf simple in $P$} if $\{ a_1, \dots,a_k\}$ is
multiplicatively independent and the subgroup $\la a_1,\dots,a_k\ra$
is pure in $P^{\times}.$ \edf

\rem In case $P=\Q$ any tuple of distinct primes is simple.\erem

\bl \lb{l1} Let $\{ a_1,\dots, a_r\}$
be a multiplicatively independent tuple  in  $P$ and  $C=P\cap \mu. $
Then the following conditions are equivalent:

(i) $\{ a_1,\dots, a_r\}$
 can be
extended to a basis of a free subgroup $A$ such that $A\cdot C=P^{\times}.$

(ii)
$\la a_1,\dots, a_r\ra $ is a direct summand of $P^{\times}.$

(iii) the image of $\{ a_1,\dots, a_r\}$ in the quotient-group
$P^{\times}/C$ can be extended to a basis of the free abelian group $P^{\times}/C.$

(iv)  $\{ a_1,\dots, a_r\}$ is simple.
\el
\pf Let $P^{\times}=A\cdot C$  and $A$ a free abelian subgroup of $P^{\times}.$ This 
induces the projection $\pr: P^{\times}\to A$ with  kernel $C.$  

Suppose now
$P^{\times}=\la a_1,\dots, a_r\ra \cdot B,$ direct. Then  $\pr$ is a
monomorphism on $\la a_1,\dots, a_r\ra $ and $\pr(A)=A$ is a direct product
of $gp(\pr(a_1),\dots, \pr(a_r))$ and $\pr(B).$ Choose now a subset $S\subs
\pr(B)$ freely generating $\pr(B)$ and $T\subs \pr\inv(S)$ such that
$\pr(T)=S$ and $\pr\inv(s)\cap T$ consists of  one element for any
$s\in S.$ Then $T$ generates a free subgroup $\la T\ra $ such that
$\la T\ra \cdot C=B$ and so $T$ completes $\{ a_1,\dots, a_r\}$ to a set
of free generators.  This proves that 
(ii) $\Rightarrow$ (i).

(i) $\Rightarrow$ (iii) by definition, and   
(iii) $\Rightarrow$ (iv) is obvious.

Also if (iii) holds, let  $U$ be  a basis of $P^{\times}/C$ extending the
image of  $\{ a_1,\dots, a_r\}$    
and $T\subs P^{\times}$ a set of representatives of elements of $U,$ containing
 $\{ a_1,\dots, a_r\}.$ It is easy to see that $T$ generates a free
group complementary to $C.$ This proves that  
(iii) $\Rightarrow$ (i).

Suppose that (iv) holds. 
Consider a basis $U$ of free generators of  $A.$ Let $u_1,\dots, u_n$ be distinct
elements of $U$ such that   $\{ a_1,\dots, a_r\}\subs  \la  u_1,\dots,
u_r\ra .$ By Corollary 28.3 in [F], under the assumptions that  
$gp\{ a_1,\dots, a_r\}$  is a pure subgroup of a finitely generated
abelian 
group, it is a direct summand of the group. I.e.
 $$\la  u_1,\dots, u_n\ra =   \la  a_1,\dots, a_r\ra \cdot B$$
for some subgroup $B$ of  $\la  u_1,\dots, u_n\ra .$ Since the latter is
free, $B$ is free too and thus  $\{ a_1,\dots, a_r\}$    can be
extended to a basis of  $\la  u_1,\dots, u_n\ra .$   Adding to the basis
$U\smin  \{ u_1,\dots, u_n\}$ we get a basis of $A.$             
 Thus 
(iv) $\Rightarrow$ (iii).          

Finally, notice that (ii) $\Rightarrow$ (iii) is obvious.
 \qed

\bl  \lb{l2} Let $\{ a_1,\dots,
a_r,a_{r+1}\}$ be simple in  $P_{r+1}=P(a_1,\dots, a_{r+1}).$ Then
$a_1,\dots, a_r$ is simple in $P_r=P(a_1,\dots, a_r).$\el
\pf Indeed, if, for $b\in \la a_1,\dots,a_r\ra,$ the equation $x^n=b$ has
a solution in $P(a_1,\dots,a_r),$ it has one in $P(a_1,\dots,a_r,
a_{r+1}),$ and hence it has a solution in the free group
$\la a_1,\dots,a_r,a_{r+1}\ra.$ It follows that the equation has a
solution in $\la a_1,\dots,a_r\ra .$\qed

\df Given a number $k>1,$ a non-zero $a\in P$ is said to be $k$-simple
if $a\notin \mu$ and for any $b\in P,$ $\epsilon\in \mu$ and an
integer $d$ $$a^d=b^{k}\cdot \epsilon\mbox{ implies } k|d.$$\edf

\rem Obviously, every simple $a\in P$ is $k$-simple. On the other hand,
e.g. in $\Q,$ $5^2$ is $3$-simple but not $2$-simple.\erem

\bl \lb{0.1-} Let $a\in P$ and $k>1,$ an integer. Then the following
three conditions are equivalent

(i) $a$ is $k$-simple;

(ii) given a  divisor $m>1$ of $k$
there is no $\alpha\in P$ and a root of unity $\epsilon$ 
such that $a=\alpha^m\epsilon;$

(iii) given a  divisor $m>1$ of $k$ the image of $a$ in the
quotient-group $P^{\times}/C$ ($C=P\cap \mu$) has no roots of power $m.$
\el
\pf (i) $\Leftrightarrow$ (iii) since the group $P^{\times}/C$ is torsion-free.

(i) $\Rightarrow$ (ii) is obvious. To prove the converse
suppose the negation of (i) holds, i.e.
 $a^d=b^k\epsilon$ for some $b\in P,$ $\epsilon \in C,$ 
$(k,d)=s<k,$ $m=k/s$  and $d$ is minimal for all choices of $b,$
$\epsilon$ and $k.$ Then by minimality $s=1$ and
$1=ku+dv$ for some integers $u$ and $v.$ Thus
$$a=a^{dv}a^{ku}=b^{kv}a^{ku}\epsilon'=(a^ub^v)^k\epsilon'$$
for some $\epsilon'\in C.$
Hence, letting $\alpha=(a^ub^v)$  we get the negation of (ii). 
\qed   

Let from now on $q$ be a prime number.

\bl  \lb{0.0}  Let $\epsilon$  be
 a root of unity of order $q.$ Suppose  $a\in P$ is $q$-simple in $P.$ Then $a$ is $q$-simple in  $P(\epsilon).$\el

\pf Suppose $$a=\alpha^q\zeta$$ for some $\alpha,\zeta\in P(\epsilon)$ and 
$\zeta^M=1$ for some integer $M.$
Then $$a^M=\alpha^{qM}.$$
The orbit of $\alpha$ under the Galois group $(P(\epsilon):P)$ consists of
$d=\deg(\alpha/P)\le (P(\epsilon):P)\le q-1$ elements of the form 
$\alpha\xi,$
for $\xi\in \mu_{qM}.$
Hence the norm has the form $$N_P(\alpha)=\alpha^d\xi'=b\in P$$
for some $\xi'\in \mu_{qM}.$
Thus $a^d=b^{q}\xi''$ for some $\xi''\in \mu_{qM}.$
 This contradicts $q$-simplicity.
\qed

\bl \lb{0.2*} If $a$ is $q$-simple in $P$ and $\imath=\sqrt{-1}\in P$ then $a$ is
$q$-simple in $P(\zeta),$ where $\zeta$ is a root of unity of order
${q^t},$  $t\ge 1.$\el
\pf By Lemma~\ref{0.0} the statement  holds for $t=1.$
Suppose it holds for $t=t_0$ and fails for $t=t_0+1.$
We may then assume that $\zeta_0,$ a primitive root of unity of order $q^{t_0},$
is in  $ P.$ 

Let $\zeta$ be a primitive root of unity of order $q^t,$
$\zeta^q=\zeta_0,$
and we may assume that $\zeta\notin P.$ 
By  [L], VI,Thm 6.2,  polynomial $x^{q}-\zeta_0$ is irreducible over $P,$ thus
$$|P(\zeta):P|=q$$
By assumptions and Lemma~\ref{0.1-}  there is
$\alpha\in P(\zeta)\smin P$ and $\epsilon\in \mu$  such that 
\be \lb{n1} a=\alpha^q\epsilon.\ee
Of course, $\epsilon\in P(\zeta).$

If $(m,q)=1$ then $a^m$ is $q$-simple in $P$ as well, so by raising the
equation (\ref{n1}) to power $m$ we may w.l.o.g. assume that the
multiplicative order
of $\epsilon$ is of the form $q^r$ for some non-negative integer $r\le
t_0+1.$
Hence $\epsilon^q\in P$ and $$\alpha^{q^{2}}=a^{q}\in P.$$

Then, by  [L], VI,Thm 6.2 again,  
$\deg(\alpha/P)$ divides $|P(\zeta):P|,$ hence
$\deg(\alpha/P)=q.$

Let $\sigma\neq 1$ be an  element of the Galois group
$(P(\zeta):P).$

Since $\alpha \in P(\zeta),$ there are unique $c_0,\dots c_{q-1}\in P$ such that
  $$\alpha=\sum_{0\le i<q-1}c_i\zeta^i.$$
and, since $\zeta^q=\zeta_0\in P,$ 
 $$\sigma(\zeta)=\zeta\cdot \epsilon, $$
for some $\epsilon,$ 
a primitive root of unity of order $q,$ which is in $P.$ 

Thus, $$\sigma(\zeta^i)=\zeta^i\epsilon^i$$ and
$$\sigma(\alpha)=\sum_{0\le i<q}c_i\epsilon^i\zeta^i.$$
 On the other hand
\be \lb{nov1}\sigma(\alpha)=\alpha\xi, \ \ \mbox{ for some } \xi, \ \ \
\xi^{q^2}=1.\ee

If $\sigma(\xi)=\xi$ then $\xi\in P$ and 
$$\sigma(\alpha)
=\sum_{0\le
i<q}c_i\xi\zeta^i.$$
Hence for all $i$ $$c_i\epsilon^i=c_i\xi,$$ which means
$c_i=0$ for all but one $i=m,$   $$\alpha=c_m\zeta^m$$ and
$$a=b^q\zeta^{mq},\mbox{ where }b\in P,$$
which contradicts the assumption of $q$-simplicity, and we are done in
this case.

So we assume that $\xi$ is a root of order $q^2,$ not in $P.$

By definitions $\sigma(\xi)=\xi^l$ for some  $1<l<q^2$ such that
$(l,q)=1.$

By induction on $k$ the condition (\ref{nov1}) extends to
$$\sigma^k(\alpha)=\xi^{1+l+\dots+l^{k-1}}.$$
Since $\sigma$ is of order $q$ we have $\sigma^q(\alpha)=\alpha$ and   
  thus \be \lb{nov2} 1+l+\dots+l^{q-1}\equiv 0(mod\ q^2).\ee
But
$$1+l+\dots+l^{q-1}=\frac{l^q-1}{l-1}$$ and hence 
$$l^q\equiv 1(mod\ q)
 \mbox{ \ \ and  \ \ }l\equiv 1(mod\ q).$$
It follows $l=q+1$ and $l^k\equiv kq +1(mod\ q^2),$ for $k=0,\dots, q-1,$ thus 
$$1+l+\dots+l^{q-1}\equiv q(1+\frac{q-1}{2})(mod\ q^2).$$
Comparing with (\ref{nov2}) we see that only $q=2$ is possible. But
in this case $\xi=\imath\in P,$  the contradiction.
\qed

From now on we assume that $\imath\in P.$\\

 Let $P_0$  be a
maximal purely transcendental extension of $\Q$ in $P,$ and
$\varphi$ the Euler function.

\bl \lb{0.4} If $a$ is $q$-simple in $P$ but not $q^w$-simple  in $P(\xi)$ for some root of
unity $\xi$ and a positive integer $w,$  then $\varphi(q^w)$
 divides $|P:P_0|.$
\el
\pf First consider the case when  $\xi$ is a primitive root of unity of
order $p,$ some $p\in \N,$ $(p,q)=1.$

Suppose \be \lb{pq1} \alpha^{q^w}=a\epsilon^m,\ee
$\epsilon$ is a primitive root of unity of order $q^w,$
 $\alpha\in P(\xi),$   and  $m\in \N.$ 

By  [L], VI,Thm 6.2,  polynomial $x^{q^w}-a$ is irreducible over $P,$ thus
$$\deg(\alpha/P)=q^w$$
and for any $\epsilon,$ a  root
of unity of order $q^w,$
there is $\sigma$ in the Galois group of the normal extension
$(P(\xi):P)$ such that $$\sigma(\alpha)=\alpha\epsilon.$$
It follows that a primitive root
of unity of order $q^w,$ denote it  $\epsilon,$ belongs to $P(\xi).$

Now we compare the degrees of some Galois extensions: 
$$|P(\epsilon\xi):P_0|=|P(\epsilon\xi):P|\cdot |P:P_0|,$$
$$
|P(\epsilon\xi):P_0|=|P(\epsilon\xi):P_0(\epsilon\xi)|\cdot|P_0(\epsilon\xi):P_0|.$$
But $|P(\epsilon\xi):P_0(\epsilon\xi)|=d_1$ is a divisor of $|P:P_0|$
([L], VI, Thm 1.12) and $$|P_0(\epsilon\xi):P_0|=\varphi(q^w)\cdot \varphi(p).$$     
Hence $$|P(\epsilon\xi):P|=|P:P_0|\inv\cdot d_1\cdot  \varphi(q^w)\cdot \varphi(p).$$  
Analogously we get
$$|P(\xi):P|=|P:P_0|\inv\cdot d_2\cdot   \varphi(p)$$ for some divisor
$d_2$ of $|P:P_0|.$

So, under the condition $\epsilon \in P(\xi),$ we obtain
$$d_1\cdot \varphi(q^w)=d_2,$$
which implies that $\varphi(q^w)$ divides $|P:P_0|.$

In the general case assume $\alpha\in P(\epsilon\xi)$ and (\ref{pq1})
holds.
Let $\beta\in P(\xi)$ satisfy the equation $$\beta^{q^u}=a\epsilon^n$$
for maximal possible integer $u$ and  some integer $n.$ Then
$\varphi(q^u)$ divides $|P:P_0|$ by the above proved, and $\beta$ is
$q$-simple in $P(\xi).$ By Lemma~\ref{0.2*} $\beta$ is $q$-simple in
$P(\xi,\epsilon),$ which implies $u\ge w$ and
$\varphi(q^w)$ divides $|P:P_0|.$   
\qed

\bc \lb{c2} If $a$ is simple in $P$ then there is a positive integer $N$
depending on $P$ such that $a^{\frac{1}{N}}\in P(\xi),$ for some root of unity $\xi,$ \
 and $a^{\frac{1}{N}}$ is simple in
$P(\xi')$ for any root of unity $\xi'$ such that $P(\xi)\subs
P(\xi').$\ec
\pf 
  Let  $N$ be the maximal positive integer with the property that
$\varphi(N)$ divides $|P:P_0|$ and $a^{\frac{1}{N}}\in P(\xi)$ for some
$\xi.$

Such an integer exists because the first part of the condition is
satisfied by at most finitely many integers.

If there is $M$ and $\xi'$ such that 
$$(a^{\frac{1}{N}})^{\frac{1}{M}}\in P(\xi'),$$
then, by  Lemma~\ref{0.4}, $\varphi(N\cdot M)$ divides $|P:P_0|.$ By the choice of
$N,$ $N\cdot M\le N,$ hence $M=1.$\qed

\bl \lb{serv} Let $A$ be a free abelian subgroup of rank $r$ of a
torsion free group $A'.$ Suppose there is a natural number $N$ such
that $a^N\in A$ for any $a\in A'.$ Then $A'$ is a 
free abelian group of rank $r.$\el
\pf By  assumptions group $A'/A$ is periodic, of a bounded exponent. 
Proposition 18.3 of [F] states under these conditions that
$A'$ is a direct product of cyclic groups, in our case all the cyclic
groups are infinite, i.e. the group is free. The rank of $A'$ must be
$r$ too, because any $b_1,\dots, b_s\in A'$ with $s>r$ must be
multiplicatively dependent, since
$b^N_1,\dots b^N_s$ are dependent elements in $A.$ \qed

\bl \lb{S} Let $a_1,\dots, a_r$ be a simple tuple in $P$ and $A$ the
subgroup of $P^{\times}$ generated by   $a_1,\dots, a_r.$ 
Let $\bar P=P(\mu),$ the extension of $P$ by all the  roots of unity, and
$A^{\#}$ be the pure hull of $A$ in ${\bar P}^{\times},$ i.e. the group of the
elements  $b\in {\bar P}^{\times}$ such that $b^n\in A$ for some $n.$ 

Then $A'=A^{\#}/A^{\#}\cap \mu$ is a free abelian group of rank $r.$ \el 
  \pf By assumptions $A$ can be naturally identified with $A/A\cap \mu$ and is a free abelian group of rank $r.$ By
  Corollary~\ref{c2}, given $b\in A^{\#}$ the least positive integer $N$
  such that $b^N\in A\cdot \mu$ is  bounded by
  $|P:P_0|.$ Thus
$A'/A$ is periodic of bounded exponent.
It follows that $A$ and $A'$ satisfy the assumptions of
  Lemma~\ref{serv}. Thus $A'$ is free abelian of rank $r.$\qed

Let from now on $\bar P=P(\mu),$ the extension of $P$ by all the  roots of unity.\\

%\df A tuple  $\{ a_1,\dots, a_r\}\subs P$ is said to be {\bf strongly
%simple} if the tuple is multiplicatively independent and $\la a_1,\dots, a_r\ra $ is pure in  the multiplicative group of field 
%$\bar P.$\\

It follows from Lemma~\ref{S} 
\bc \lb{c4} For any multiplicatively
independent  tuple  $\{ a_1,\dots, a_r\}\subs P^{\times}$ there is a
tuple  $$\{ a'_1,\dots, a'_r\}\subs {\bar P}^{\times}\cap  a_1^{\Q}\cdot\dots\cdot a_r^{\Q}    $$ such that
 $\{ a'_1,\dots, a'_r\}$ is simple in the field 
$\bar P.$ \ec

The following statements \ref{root} - \ref{l3} make  sense in the case
$n\neq 0.$

\bl \lb{root} Given a field $K$ containing all roots of unity, suppose  $\{ a_1,\dots,a_r\}\subs K^{\times}$ is simple in
$K$  and $b_1,\dots, b_r$ in an extension of
$K$ generate a subgroup $\la b_1,\dots, b_r\ra $ containing  $\la a_1,\dots, a_r\ra .$

Then
$\{ b_1,\dots, b_r\}$ is simple in 
$K(b_1,\dots, b_r).$\el
\pf It follows from assumptions that $(\la b_1,\dots,
b_r\ra :\la a_1,\dots, a_r\ra )$ is finite, hence there is an $m$ such that
$b_1^m,\dots, b_r^m\in \la a_1,\dots, a_r\ra .$

Claim. Suppose $b\in \la b_1,\dots, b_r\ra $ has a root 
$\beta\in K(b_1,\dots, b_r)$ of power $l.$ 
Then $b$ has a root of power $l$ in $\la b_1,\dots,
b_r\ra .$

Indeed,  there is a positive integer $M\le ml$ such that $\beta^{M}\in
K$ and
$b_i^M\in K,$ all $i.$
By [L], VI, Theorem 8.1 (Kummer's Theory)
$$ (K(\beta,b_1,\dots, b_r): K)= 
(\la \beta,b_1,\dots, b_r\ra \cdot K^{M}:K^{M}),$$
$K^{M}$ is the $M$-powers subgroup of $K^{\times}.$ 

On the other hand $K(\beta,b_1,\dots, b_r)=K(b_1,\dots,b_r)$ and 

$$(K(b_1,\dots, b_r): K)= (\la b_1,\dots, b_r\ra \cdot K^{M}:
K^{M}).$$
It follows $\beta\in \la b_1,\dots, b_r\ra \cdot K^{M}.$ 
Hence $b=\beta^l=c^l\cdot p^l$ for some $c\in \la b_1,\dots,
b_r\ra $ and $p\in K^M.$ But $b^m,$ $c^m$ are in 
$\la a_1,\dots, a_r\ra $ and $b^m\cdot c^{-lm}\in \la a_1,\dots, a_r\ra $ 
has  a root $p$ of power $lm$ in $K.$
Since $\{ a_1,\dots, a_r\}$ is simple in $K,$ we get that $p\in
\la a_1,\dots, a_r\ra $ and thus 
$\beta$ is equal, up to a root of unity, to $c\cdot p\in
\la b_1,\dots, b_r\ra .$ This proves the claim.

The Lemma now follows directly from the claim.   \qed

\bl \lb{root2} Suppose  $\{ a_1,\dots, a_r\}\subs  P^{\times}$ are multiplicatively independent
over $L_1^{\times}\cdot \dots\cdot L_n^{\times}.$ Then
there are $a'_1,\dots, a'_r\in P(L_1,\dots, L_n)$ which are free
generators of the group $a_1^{\Q}\cdot\dots\cdot a_r^{\Q}\cap 
P(a'_1,\dots, a'_r,L_i)$ for any $i\in \{ 1,\dots,n\}.$
\el
\pf Suppose, for $0\le k<n,$ we found $a'_1,\dots, a'_r\subs P(L_1,\dots, L_n)$ which are free
generators of the group $a_1^{\Q}\cdot\dots\cdot a_r^{\Q}\cap 
P(a'_1,\dots, a'_r,L_j)$ for any $j\in \{ 1,\dots,k\}.$
Since
$a_1^{\Q}\cdot\dots\cdot a_r^{\Q}=
{a'_1}^{\Q}\cdot\dots\cdot {a'_r}^{\Q},$ we may assume $a'_i=a_i\in P,$
$i=1,\dots, r.$
By Lemma~\ref{l0},  
$a_1^{\Q}\cdot\dots\cdot a_r^{\Q}\cap 
P(a_1,\dots, a_r,L_{k+1})$ is free, so let 
$a'_1,\dots, a'_r$ be its free
generators.

By Lemma~\ref{root} $\{ a'_1,\dots, a'_r\}$ is simple in 
$P(a'_1,\dots, a'_r,L_j)$ for $j\le k$ and the same holds by the choice
for $j=k+1.$ The lemma follows now by induction on $k.$
\qed

\bp \lb{ssL} Given $\{ a_1,\dots,
a_r\} \subs \bar P$   multiplicatively
independent over $\bar P\cap L_1^{\times}\cdot \dots \cdot L^{\times}_n,$  
there are 
$\{ a'_1,\dots,
a'_r\}\subs \bar P(L_1,\dots,L_n)$  simple in $\bar P(L_1,\dots,L_n),$
such that $a_1,\dots,a_r\in \la a'_1,\dots, a'_r\ra .$
\ep
\pf  By Lemma~\ref{root2} we can find $\{ a'_1,\dots,a'_n\}$  simple
in $\bar P(a'_1,\dots, a'_r,L_j)$ for each $j\in \{ 1,\dots, n\}.$
We may assume $a'_i=a_i\in \bar P,$ all $i,$ and are going to prove that 
$\{ a_1,\dots, a_r\}$  is simple in
$P(L_i)$ for each $i\in \{ 1,\dots, n\}.$
 
Since $L_1,\dots, L_n$ are countable fields, we can represent
$$P(L_1,\dots,L_n)=\bar P(L_1,\dots, L_n)=\bigcup_{i\in \N}P^{(i)},$$
where $P^{(0)}=\bar P$ and 
$P^{(i+1)}=P^{(i)}(\lambda_i)$ for some $\lambda_i\in L_1\cup \dots \cup
L_n.$ Moreover,  letting $L_j^{(i)}=L_j\cap P^{(i)},$ we may assume that
either  $L_j^{(i+1)}=L_j^{(i)},$ or $L_j^{(i+1)}=L_j^{(i)}(\lambda_i)$ and
$(L_j^{(i+1)}:L_j^{(i)})$ is a
normal extension with simple Galois group, that is with no intermediate normal extensions. 
In the second case, since $\lambda_i\in \acl(L_j^{(i)})$ and 
$L_j^{(i)}$ is algebraically closed in $P^{(i)},$ by Lemma~4.10 of [L],
VIII, we have linear disjointness and an isomorphism of Galois groups
$(L_j^{(i+1)}:L_j^{(i)})$ and $(P^{(i+1)}:P^{(i)}).$

To prove the proposition it is enough to check that if an element
$a\in \bar P$ is simple in
$P^{(i)}(L_j),$ all $j,$ it
is simple in $P^{(i+1)}(L_j)$ all $j.$ 

Suppose towards a contradiction that
$a$ is not simple in  $P^{(i+1)}(L_j).$  That is there is a root $\alpha$ of $a$ of some
order $m>1$ in $P^{(i+1)}(L_j).$ Choosing $m$
minimal, we have that $\alpha$ generates a cyclic (hence normal) extension of
$P^{(i)}$ of order $m.$ Since $(P^{(i+1)}:P^{(i)})$ has no intermediate normal
extensions, it must be cyclic generated by $\alpha,$ so 
$(L_j^{(i+1)}:L_j^{(i)})$ is cyclic of order $m$ as well. Hence, we
may assume $\lambda_i=\lambda$ is a root of order $m$ of an element
$b\in L_j^{(i)}.$ Consider a unique representation
$$\lambda=p_0+p_1\alpha+\dots+p_{m-1}\alpha^{m-1},$$
with $p_0,\dots,p_{m-1}\in P^{(i)}.$

Let $\sigma$ be an automorphism of $(P^{(i+1)}:P^{(i)})$ which sends
$\alpha$ to $\alpha\xi,$ for $\xi$ a primitive root of $1$ of order
$m.$
Then 
$$\sigma(\lambda)=p_0+p_1\alpha\xi+\dots+p_{m-1}\alpha^{m-1}\xi^{m-1}$$
and at the same time
$$\sigma(\lambda)=\lambda\xi^k=p_0\xi^k+p_1\alpha\xi^k+\dots+p_{m-1}\alpha^{m-1}\xi^k
$$ for some $k.$ 

Comparing the two expressions we get that all but one $p_i$'s is zero
and $\lambda=p_k\alpha^k,$ in fact $k$ coprime with $m.$ Thus 
$\alpha=p_k\inv\lambda^{m'}$
for some $m'$ and so $\alpha\in P^{(i)}\cdot L_j,$ proving that $a$ is not simple in $P^{(i)}(L_j).$  
The contradiction.\qed

\bl \lb{l3} Let $R\subs \C$ be a field containing a primitive root of
unity of order $n,$ 
$\{ a_1,\dots,
a_r\} \subs R$ be simple in $R$ and $\alpha_1,
\dots,\alpha_r\in \C,$ 
$\alpha_i^n=a_i$ for $1\le i \le r.$  Then the Galois group
$(R(\alpha_1,\dots, \alpha_{r}): R)$ is isomorphic to $\Z_n^r,$ 
the $r$th
Cartesian power of
cyclic group of order $n.$ That is any other collection
$(\alpha'_1,\dots,\alpha'_r)$ of roots of $(a_1,\dots,a_r)$ of order
$n$ is conjugated to $(\alpha_1,\dots,\alpha_r)$ by an automorphism over $R.$\el 
\pf 
Let ${R^{\times}}^n$ be the $n$-powers subgroup of $R^{\times}.$ Since
$\{ a_1,\dots,a_r\}$ is simple, we have  the group isomomorphism
$$\la a_1,\dots,a_r\ra /\la a_1,\dots,a_r\ra \cap {R^{\times}}^n\cong\Z_n^r.$$
On the other hand, by Kummer's theory (Theorem 8.1 of [L], Ch.VI)  
$$(R(\alpha_1,\dots, \alpha_{r}): R)\cong
\la a_1,\dots,a_r\ra /\la a_1,\dots,a_r\ra \cap {R^{\times}}^n.$$ 
\qed

\pf of the theorem. 

Let $C=\mu\cdot L_1^{\times}\cdot\dots\cdot L_n.$
We may assume that $\{ a_1,\dots, a_r\}$ is
multiplicatively independent over $C.$ Let $$A= \mu\cdot  \la a_1,\dots,
a_r,b_1,\dots, b_l\ra \cap P(b_1,\dots,b_l).$$ Since $A/A\cap C$ is a free group 
of rank $r+l,$ by Proposition~\ref{ssL} there is  \\
$\{ a'_1,\dots,a'_r,b'_1,\dots, b'_l \}\subs a_1^{\Q}\cdot \dots \cdot 
a_r^{\Q}\cdot b_1^{\Q}\cdot\dots\cdot b_l^{\Q}$  simple
in $\bar P(a'_1,\dots,a'_r,b'_1,\dots, b'_l).$ We preserve this property
with any choice of free generators of the group
$\la  a'_1,\dots, a'_r,b'_1,\dots,b'_l\ra ,$ in particular we may assume  
$$\{ a'_1,\dots, a'_r\}\subs   a_1^{\Q}\cdot\dots\cdot a_r^{\Q}.$$     
 By Lemma~\ref{l3} $(b'_1,\dots, b'_l)$
determines the type of $({b'_1}^{\Q},\dots, {b'_l}^{\Q})$ over \\
$\bar P({a'_1}^{\Q},\dots, {a'_r}^{\Q},L_1,\dots,L_n).$

Obviously, $b'_1,\dots, b'_l$ are in the subgroup generated by
$b_1^{\frac{1}{m}},\dots, b_l^{\frac{1}{m}},$  $\mu$ and 
$a_1^{\Q}\cdot\dots\cdot a_r^{\Q},$
for some integer  $m,$ and thus $(b_1^{\frac{1}{m}},\dots, b_l^{\frac{1}{m}})$
determines the type of $({b_1}^{\Q},\dots, {b_l}^{\Q})$ over  $\bar P({a'_1}^{\Q},\dots, {a'_r}^{\Q},L_1,\dots,L_n).$
    \qed

\ssn{ Theorem~\ref{c1}}

As was mentioned in the first section of the paper, the proof of
Theorem~\ref{c1} from Theorem~\ref{t0} is based on a rather general
model theoretic construction.

First, let us represent a sequence of the form (\ref{e4}) as a
one sorted structure $\HH.$ The domain of $\HH$ will be $H,$ the 
only basic operation is $+,$ the group operation, and the basic
relations on $H$ are: a
binary
equivalence relation $E,$ with interpretation
$$E(h_1, h_2) \mbox{ iff } \ex(h_1)=\ex(h_2),$$
     and a ternary  relation $S$  with the interpretation
$$S(h_1,h_2,h_3)\mbox{ iff } \ex(h_1)+\ex(h_2)=\ex(h_3).$$

So, each equivalence class $hE$ with a representative $h\in H$
corresponds to a non-zero element $\ex(h)$ of $F.$ The
multiplication in $F$ corresponds to the group operation $+$ on $H,$
and we also have $S$ to 
speak about addition in $F.$ In particular, the equivalence class corresponding to
the unit $1$ of $F,$ which is just the kernel of $\ex,$  is definable in $\HH.$  

The appropriate formal 
logical language to consider the structures  is $L_{\omega_1,\omega},$
the language 
with countable conjunctions and finite number of variables studied in [K].

\bl \lb{l5} 
 There is an $L_{\omega_1,\omega}$-sentence $\Sigma$  such that any
model $\HH$ of $\Sigma$ represents a sequence (\ref{e4}) with some
algebraically closed field $F$ and conversely, any $\HH$ corresponding
to a sequence of the form (\ref{e4}) is a model of $\Sigma.$
\el
\pf $\Sigma$ should say that: 

(i) $(H, +)$ is a divisible torsion free
 abelian group; 

(ii) $H/E,$ with regards to the operations coming from $+$ and $S,$ can be
 identified as $F\smin \{ 0\}$ for some algebraically closed field $F$ of
 characteristic zero; 

and 

$${\rm (iii)}\mbox{ \ \ \ \ \ \ \ \ \ \ \ \ \ \ \ \ \ \ \ \ }\exists
x_0\in \Ker, \forall
x\in\Ker \ \ \
( \bigvee_{z\in \Z} x=z x_0),\mbox{ \ \ \ \
\ \ \ \ \ \ \ \ \ \ \ \ \ \ \ \ \ \ }$$
where $\Ker $ stands for the kernel of the homomorphism $\ex.$
\qed

We also can define a closure operator $\cl$ on the domain $H$ of a
model $\HH$ of $\Sigma$ by letting for $X\subs H$
$$\cl(X)=\ex\inv(\acl(\ex(X))),$$
where $\acl$ is the algebraic closure operator in the sense of the
field structure on $F=\ex(H)\cup \{ 0\}.$

\bl \lb{I}   For any model $\HH$ of $\Sigma$ and  $X\subs H:$

(i) $\cl(X)$ is countable for $X$ finite;\\ 

(ii) $$\cl(Y)=\bigcup_{X\subs Y, \ \ X\mbox{ finite}} \cl(X);$$

(iii) $X\to \cl(X)$ is a monotone idempotent operator;

(iv) $\cl$ satisfies the exchange principle:
$$ z\in \cl(X\cup \{ y\})\smin \cl(X) \Rightarrow y\in \cl(X\cup\{ z\});$$

(v) $\cl(X)$ with the induced relations is a model of $\Sigma.$
\el
\pf Follows immediately from the properties of $\acl.$\qed 

\df Let $\HH,\HH'$ be models of $\Sigma$ and $G$ their common  subset.
A (partial) mapping  $\varphi: \HH\to \HH'$ is called a
{\bf $G$-monomorphism}, if it preserves quantifier-free formulas with
parameters from $G,$ that is for any such formula
$\Phi(v_1,\dots,v_n)$ and elements $h_1,\dots,h_n\in H$
$$\HH\models \Phi(h_1,\dots,h_n) \mbox{ iff } \HH'\models \Phi(\varphi(h_1),\dots,\varphi(h_n)).$$\edf

\rem The $G$-monomorphism type of a linearly independent tuple 
$(x_1,\dots,x_l)$ in $H$ is
determined by the algebraic type of
$(\ex(q_1x_1),\dots,\ex(q_lx_l))$ for all  rational numbers
$q_1,\dots,q_l.$ In other words, letting 
$y_i^{q}=\ex(qx_i),$
we want to know the field theoretic isomorphism type of
$(y_1^{\Q},\dots, y_l^{\Q}).$\erem

\bl \lb{II} Models of $\Sigma$  are {\bf $\omega$-homogeneous over a
model}. That is, given $\HH$ and $\HH'$ models of $\Sigma$ and  a
common submodel $G\subs \HH,$ $G\subs \HH'$ the following holds:

(i) Suppose $X_0\subs H,$ $X'_0\subs H'$  are finite subsets of models
$\HH$ and $\HH'$ correspondingly, and there is a $G$-monomorphism
$\varphi_0: X_0\to X'_0.$  Suppose also 
$X\subs H$ and $X'\subs H'$ are  $\cl$-independent  over
$X_0\cup G$ and $X'_0\cup G,$ 
 correspondingly. Then any bijection  $\varphi:
X_0X\to X'_0X'$ extending $\varphi_0$  is a  $G$-monomorphism.\\

(ii) If a partial $\varphi: \HH
\to \HH'$ is  a $G$-monomorphism, $\Dom \varphi = X,$ with $X$ finite, then
for any $y\in \HH$ there is a $G$-monomorphism  $\varphi'$ extending $\varphi$ with
$\Dom \varphi' = X\cup \{ y\}.$  \\

(iii) if $\varphi: X\cup \{y\}\to X'\cup \{ y'\}$ is a  $G$-monomorphism,
then  $$y\in \cl(X) \mbox{ iff }y'\in \cl(X').$$ 
\el
\pf (i) is obvious if one remembers that the $\cl$-independence of $X$
over $X_0$ means that $\ex(X)$ is algebraically independent over
$\ex(X_0)$ in the field $F.$  

(ii) follows directly from Theorem~\ref{t0} (see also
Corollary~\ref{c0}), when one takes $\ex(G)=L=L_1=\dots=L_n,$
algebraically closed subfield of $F,$  $X$ to be $\{ a_1,\dots, a_r\}$
and $y=b_1,$ $l=1.$ 

(iii) is obvious.\qed

\bl \lb{III} Given a countable submodels $G_1,\dots, G_n\subs \HH$ and  
$h_1,\dots, h_l\in \cl(G_1\cup\dots \cup G_n),$ the locus of $(h_1,\dots,h_l)$
over $G_1\cup\dots \cup G_n$ is finitely determined, i.e. there is a finite subset
$A\subs G_1\cup\dots \cup G_n,$ such that any $\varphi: \{ h_1,\dots, h_l\}\to \HH$ which is an
$A$-monomorphism is also a $(G_1\cup\dots \cup G_n)$-monomorphism.  
\el
\pf Let $L_i=\ex(G_i),$ $i=1,\dots,n,$
$b_j^{q}=\ex(qh_j),$ for $j=1,\dots,l,$ $q\in \Q.$ We may
assume that $h_1,\dots, h_l$ are $\Q$-linearly independent over the vector
subspace $G_1+\dots+G_n,$ which implies the multiplicative
independence of $b_1,\dots,b_l$ over $L_1^{\times}\cdot\dots\cdot L_n^{\times}.$

Apply Theorem~\ref{t0} with the
above notation assuming $r=0.$ Since the type of $(h_1,\dots,h_l)$
over $G_1\cup\dots \cup G_n$ is determined by the field-theoretic type of 
$(b_1^{\frac{1}{m}},\dots,b_l^{\frac{1}{m}}),$ for some $m,$ only
finitely many parameters from $G_1\cup\dots \cup G_n$ are needed to
fix the type. 
\qed  

In [Z2] we call an $L_{\omega_1, \omega}$-definable class  {\bf
quasi-minimal excellent} if it satisfies the statements of
Lemmas \ref{I}-\ref{III}.\\

\pf of Theorem~\ref{c1}. The main Theorem 2 of [Z2] immediately implies that, if the class of models of a
$\Sigma$ is quasi-minimal excellent, then, given an
uncountable cardinality, a
model of $\Sigma$ of this cardinality is unique,   up to
isomorphism. This yields the proof of Theorem~\ref{c1}.     \qed

\rem \lb{r3.7} The unique model of $\Sigma$ of cardinality $\omega_1$ is not
 homogeneous. So we can not
apply Keisler's theory of categoricity and stability (see [K]) to this
 sentence. In fact, $\Sigma$ provides a natural example for
the negative answers to Open Question in [K], pp. 100-101. Earlier an
 artificial 
 counterexample to the questions was published in [M].\erem

 Indeed,
consider a transcendental $a\in F^{\times}$ and $h\in H$ such that $\ex(h)=a.$ Let for every
$n\in \N$ $$a_n=\ex(\frac{1}{n}h)+1\mbox{ and } \ex(h_n)=a_n.$$
Let $X=\{ h_n: \ n\in \N\},$ \ $X_i=\{ h_n: \ n\le i\}$ \ 
and $$p=\tp(h/X).$$
Notice that $a,a_1,\dots,a_n$ are multiplicatively independent over
$\Q^{\times},$ since $a$ is transcendental. 
Then any subtype $p_i=\tp(h/X_i)$ is, by Theorem~\ref{t0}, atomic and
actually defined by the minimal polynomial for 
$a^{\frac{1}{m}}=\ex(\frac{h}{m}),$ for some $m,$ over $\Q(\ex(span_{\Q} X_i))$ which is also
an $L_{\omega_1,\omega}$-complete formula. This  also implies that any two
roots of $a^{\frac{1}{m}}$ of any power $k>0$ are indiscernible over
$X_i.$ But they are discernible over $X_{mk},$ being elements of  
$\Q(\ex({span}_{\Q} X_{mk})).$
Hence  $p_i$ is not complete over $X,$ hence $p$ is not atomic. It
follows that, given any countable fragment $L'$ of $L_{\omega_1,\omega},$  there is a countable model $\HH_0$ of $\Sigma$ containing
an $L'$-equivalent copy $X'$ of $X$ and omitting type $p',$ which
corresponds to  $p.$ 
By $\omega_1$-categoricity $\HH_0$ is embeddable in $\HH,$ and 
$p'$ is omitted in $\HH$ as well, since any realisation $h'$ of $p'$ satisfies
$\ex(h')\in \acl(\ex(X))\subs \HH_0$ and hence $h'\in \HH_0.$   
 \\ 

The last observation contrasts with the case of a cover with {\em
compact kernel}, that is the case when the kernel of $\ex$ is assumed
to be $\proZ,$ the closure of $\Z$ in the profinite topology,
corresponding to the sequence of the form
\be \lb{ep1} 0\to \proZ\to^{i_H} H\to^{\ex_H} F^{\times}\to 1. \ee

\bp \lb{pl} If in an exact sequence of groups
\be \lb{ep2} 0\to \proZ\to^{i_G} G\to^{\ex_G} F^{\times}\to 1 \ee 
 $G$ is a divisible torsion-free abelian group, then there is a group isomorphism 
$$\sigma: H\to G$$ such that $\ex_H=\ex_G\circ \sigma.$\ep 
\pf Let $h\in H$ and $g\in G$ satisfy $\ex_H(h)=\ex_G(g)=a$ for some
$a\in F^{\times}.$ \\

{\bf Claim.} There is a unique $\nu\in \Ker \ex_G$ such that $$\ex_H(\frac{h}{n})= 
\ex_G(\frac{g+\nu}{n})\mbox{ for all } n\ge 1.$$
Indeed, let $$b_n=\ex_H(\frac{h}{n})\ex_G(\frac{g}{n})\inv\mbox{ and }
\beta_n\in G\mbox{ such that } \ex_G(\frac{\beta_n}{n})=b_n,\mbox{ for
each }n\ge 1.$$
Obviously, $$\beta_n\in \Ker \ex_G\mbox{ and }
\frac{\beta_{nm}-\beta_n}{n}\in \Ker \ex_G.$$
We may identify the additive group $\Ker \ex_G$ with the group $\proZ$
and 
look for $\nu$ as a solution of the system of equations mod $n$ for
$z\in  \proZ$  

\be \lb{proZ}z\equiv \beta_n\ \mod n\proZ, \ \ \ \ n\in \N.  \ee
The system is finitely satisfiable, for to find a
solution for  $n_1,\dots,n_k,$ it is enough to solve one equation 
$$z\equiv \beta_m\ \mod m\proZ$$
for $m=n_1\cdot\dots\cdot n_k.$

The defining property of $\proZ$ is that any finitely satisfiable system of the form
(\ref{proZ}) has a unique solution (see [F], Ch.VII). This proves the claim.\\

By the Claim to every $h\in H$ we can assign unique $g$ with the
property $$\ex_H(\frac{h}{n})= 
\ex_G(\frac{g}{n})\mbox{ for all } n\ge 1.$$
Letting $\sigma(h)$ we have   $\ex_H=\ex_G\circ \sigma$ and also, by
uniqueness, $\sigma(h_1+h_2)=\sigma(h_1)+\sigma(h_2).$ \qed
{\bf References}\\

[F] L.Fuchs, {\bf Infinite Abelian Groups}, vol.1, Academic Press, New
York and London, 1970\\

[H] R.Hartshorne, {\bf Algebraic Geometry}, Springer-Verlag, Berlin,
Heidelberg, New York, 1977\\

[He] H.Hasse, {\bf Number Theory},Grundlehren der mathematischen 
Wissenschafte; 229, Springer-Verlag, Berlin, 1980\\

[L]  S.Lang, {\bf Algebra}, Addison-Wesley Publishing Co., Reading,
Mass., 1965\\

[M] L.Marcus, {\em A prime minimal model with an infinite set of
indiscernibles}, Isr.J.Math., 1972, 180-183\\

[K] H.J. Keisler, {\bf Model theory for infinitary logic. Logic
with countable conjunctions and finite quantifiers}. 
Studies in Logic and the Foundations of Mathematics, Vol. 62. 
North-Holland Publishing Co., Amsterdam-London, 1971\\

[Sh] S.Shelah, {\em Classification theory for non-elementary classes.I:
the number of uncountable models of 
$\psi \in L_{\omega_,\omega}(Q).$}, Parts A,B,  Isr.J.Math., v.46,
(1983), 212-240, 241-273.\\

[Z1] B.Zilber, {\em   Analytic and pseudo-analytic structures}, To
appear in {\bf Proceedings of European Logic Colloquium. Paris 2000}  \\ 

[Z2] --------, {\em A categoricity theorem for quasi-minimal excellent classes,
} 
Available on www.maths.ox.ac.uk/\~{ }zilber\\ 

[Z3] --------, {\em Pseudo-exponentiation on algebraically closed
fields of characteristic zero,} Submitted\\

[Z4] --------, {\em Model theory, geometry and arithmetic of
semi-Abelian varieties} In preparation.
\end{document}